\documentclass[dvips,11pt]{amsart}
\usepackage{psfrag}

\newcommand{\dsp}{\displaystyle}
\newcommand{\qq}{{\bf q}}
\newcommand{\ff}{{\bf f}}

\newcommand{\uu}{{\bf u}}
\newcommand{\vv}{{\bf v}}
\newcommand{\bb}{{\bf b}}
\newcommand{\IR}{\mathbb{R}}
\newcommand{\IN}{\mathbb{N}}
\newcommand{\uueps}{{\bf u}^\varepsilon}
\newcommand{\qqeps}{{\bf q}^\varepsilon}
\newcommand{\ueps}{{u}^\varepsilon}
\newcommand{\qeps}{{q}^\varepsilon}
\newcommand{\eps}{\varepsilon}
\def\virg {\;,\;\;}
\def \PP {{\mathcal P}}
\def \tPP {{\mathcal{ \tilde P}}}
\def \QQ {{\mathcal Q}}
\def \tQQ {{\mathcal{ \tilde Q}}}
\def \ONE {\mathbf{1}}

\def \vseq {\vspace{7pt}}
\def \qqinit {\qq_0}

\newcommand{\Linf}{{L^\infty(I)}}

\newcommand{\Hunzero}{{H^1_0(I)}}
\newcommand{\LunLun}{{L^1(]0,T[\times I)}}
\newcommand{\LdeuxLdeux}{{L^2(]0,T[\times I)}}

\newcommand{\LinfLun}{{L^\infty(]0,T[,L^1(I))}}
\newcommand{\LinfLdeux}{{L^\infty(]0,T[,L^2(I))}}
\newcommand{\LinfLinf}{{L^\infty(]0,T[\times I)}}
\newcommand{\HunHun}{{H^1(]0,T[\times I)}}

\newcommand{\LinfHun}{{L^\infty(]0,T[,H^1_0(I))}}

\newcommand{\wsto}{\stackrel{\star}{\relbar\joinrel\rightharpoonup}}

\theoremstyle{plain}
\newtheorem{theorem}{Theorem}[section]
\newtheorem{lemma}[theorem]{Lemma}
\newtheorem{proposition}[theorem]{Proposition}
\theoremstyle{definition}
\newtheorem{definition}[theorem]{Definition}
\newtheorem{remark}[theorem]{Remark}

\begin{document}

\markboth{A. Lefebvre, B. Maury}{Micro-Macro Modelling of an Array of Spheres Interacting Through Lubrication Forces}

%
%

\title[MICRO-MACRO MODELLING]{MICRO-MACRO MODELLING OF AN ARRAY OF SPHERES\\ INTERACTING THROUGH LUBRICATION FORCES}

\author{A. LEFEBVRE-LEPOT}


\author{B. MAURY}


\thanks{\begin{tabular}{rl}{\it Address:}
& A. Lefebvre-Lepot, CMAP, CNRS et Ecole Polytechnique, route de Saclay, 91128 Palaiseau Cedex, France
\\&B. Maury, Laboratoire de Math\'ematiques d'Orsay, Universit\'e Paris-Sud 11, 91405 Orsay Cedex, France
\\&Aline.Lefebvre@cmap.polytechnique.fr, Bertrand.Maury@math.u-psud.fr\end{tabular}
}

\maketitle


\begin{abstract}
We consider here a discrete system  of spheres interacting through a
lubrication force. This force is dissipative, and singular near
contact: it behaves like the reciprocal of the interparticle distance.  
We propose a macroscopic constitutive equation which is  built as the natural continuous counterpart of this microscopic lubrication model. This model, which is of the newtonian type, relies on an extensional viscosity, which is proportional to the reciprocal of the local fluid fraction. We then establish  the convergence in a weak sense of solutions to the discrete problem towards the solution to the partial differential equation which we identified as the macroscopic constitutive equation.

\end{abstract}

\begin{center}
\begin{minipage}{0.85\textwidth}
{\small
\vspace{6pt}

}
\end{minipage}
\end{center}

\section{Introduction}
We are interested in the macroscopic behaviour of highly dense suspensions of rigid spheres. 
The case of dilute suspensions has motivated  a huge  litterature since 1906, where the first derivation of the first order expansion of the equivalent viscosity with respect to the solid fraction has been proposed by Einstein (see~\cite{einstein},  or~\cite[\S  22]{landau},  where Einstein's approach is detailed).
  This approach was extended to semi-dilute
suspensions (see~\cite{batchelor}) which leads to second order asymptotic expansions of the apparent viscosity with respect to the solid fraction. 
For intermediate volume fractions (say, between 5 and 30 \%), particles interact in a complex way, so that direct numerical simulation is usually considered as the only way to account for the complex phenomena which are likely to occur at the microscopic level (see e.g.~\cite{glo}, or ~\cite{mauryapp} for an illustration of how direct simulation can be used to investigate the apparent viscosity of a suspension of particles at intermediate volume fraction).

For highly packed suspensions in a viscous fluid, interparticle distances  tend to approach zero, so that lubrication forces between particles in quasi-contact  become predominant. It is then natural to consider
 the suspension as a collection of spheres, each of which interacting with its close neighbours only, according to some model accounting for the lubrication forces.

 %
 %
%
%
In Ref.~\cite{snabre}, a first attempt was proposed to investigate the
behaviour of the apparent shear viscosity of a suspension in the
neighbourhood of the maximal packing solid fraction $\Phi_{\rm
  max}$. A model is proposed, which gives a shear viscosity which is
of the order $(1-\Phi/\Phi_{\rm max})^{-2}$ where $\Phi$ is the solid fraction.
In Ref.~\cite{berlyand}, the authors investigate the asymptotic behaviour,
as $\varepsilon$ goes to $0$,  of a set of particles under the
assumption  that distances between neighbouring particles are subject
to behave like $\varepsilon$. In this framework, the authors establish
that the apparent shear viscosity behaves like
$1/\varepsilon^{3/2}$. This approach extents a previous work~\cite{frankel} where periodic arrays of spheres are considered. In
this context, the elongational viscosity can be shown to behave like
$1/d$ where $d$ is the constant distance between neighbouring spheres.

The approach we propose here is based on a simpler model from the geometric standpoint, as the spheres
are supposed to be aligned. On the other hand it generalizes these works in the
sense that no assumption is made on the distances:
the macroscopic behaviour depends on the local solid fraction only. 
Contacts between neighboring particles are even  allowed, and  a special attention has been paid to the way we express the continuous model so that macroscopic clusters  can be taken into account (the local viscosity within a cluster is infinite). 

{ This work can be related to classical results in homogeneization of elliptic operators (see e.g. \cite{murat,spagnolo} and references therein). It consists in studying the asymptotic behaviour of a sequence of solutions $u_\varepsilon$ to  elliptic problems:
$$
-\hbox{div}(a_\varepsilon (x) \nabla u_\varepsilon)=f \hbox{ in } \Omega
$$
$$
u_\varepsilon=0 \hbox{ on } \partial\Omega
$$
where $\Omega$ is a bounded set of $\IR^d$ and $a_\varepsilon$ a sequence of equi-coercive matrix-valued functions. Most of the results concerning the asymptotic behaviour of such systems suppose that the sequence $a_\varepsilon$ is bounded in $L^\infty$. In the two dimensional case, some results were obtained in~\cite{briane1,briane2} under more general assymptions. However, $a_\varepsilon$ still have to belong to $L^\infty$. In dimension one, when the diffusion coefficient $a_\varepsilon$ belongs, together with  its reciprocal,  to $L^\infty(\Omega)$, the closure of the set defined by such systems has been studied in~\cite{alibert}. In our model, such results do not hold anymore. Indeed, since we allow contacts between particles, one of the difficulties of this work is that $a_\varepsilon$ can take infinite values on non negligible subsets of $\Omega$. Note that another main difference between our work and classical homogenization results stands in the fact that we study non stationnary problems: the density evolves according to a transport equation.} 
 
We prove in Section~\ref{sec:micmacstat}  that the limit elongational viscosity behaves singularly with respects to the vanishing fluid fraction.
This approach leads to an  equation of the elliptic type
\begin{equation}
\label{macro_ss_masse}
- \partial_x \left( \frac{1}{1-\rho} \partial_x u\right) = \rho f.
\end{equation}

where $u$ is the velocity field, $\rho$ the solid fraction (which is $1$ when all particles are in contact), and $f$ an external body force.
This result is extended to moving particle collections in Section~\ref{sec:micmacnonstat}. In this latter situation, we establish convergence to non-stationnary velocity and density fields, such that instantaneous balance  of forces reads as previously, and solid phase motion is described by the standard transport equation
\begin{equation}
\partial _t \rho  + \partial _x (\rho u) = 0.
\end{equation}
 

\section{Discrete Model}

\begin{figure}[th]
\psfragscanon
\psfrag{q1}[l]{$q_1$}
\psfrag{q2}[l]{$q_2$}
\psfrag{D}[l]{$d$}
\centerline{\resizebox{0.50\textwidth}{!}{
\includegraphics{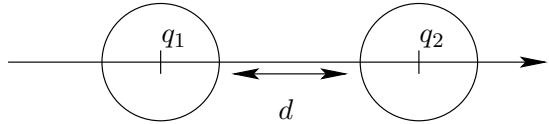}
}}
\vspace*{8pt}
\caption{Lubrication model}
\label{fig:two_spheres}
\end{figure}

\begin{figure}[th]
\psfragscanon
\psfrag{q0}[l]{$q_0$}
\psfrag{q1}[l]{$q_1$}
\psfrag{qim1}[l]{$q_{i-1}$}
\psfrag{qi}[l]{$q_i$}
\psfrag{qNm1}[l]{$q_{N-1}$}
\psfrag{qN}[l]{$q_N$}
\psfrag{d1}[l]{$d_1$}
\psfrag{di}[l]{$d_i$}
\psfrag{dNm1}[l]{$d_{N-1}$}
\psfrag{eps}[l]{$\varepsilon$}
\centerline{\resizebox{0.85\textwidth}{!}{
{\includegraphics{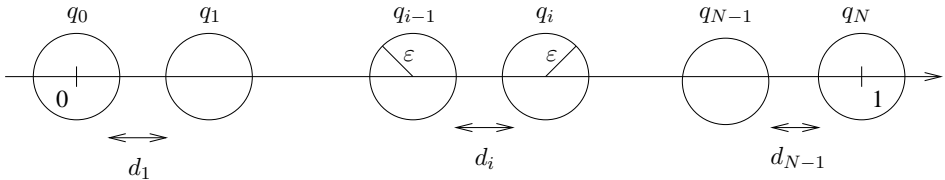}}
}}
\vspace*{8pt}
\caption{Geometry}
\label{fig:geom}
\end{figure}

Consider (see Fig.~\ref{fig:two_spheres}) two rigid spheres imbedded in a viscous fluid, subject to move horizontally.
Denoting by $q_1$ and  $q_2$ the abscisses of their centers, by $u_1$ and
$u_2$ their instantaneous velocities and by $d$ the border-to-border
distance, the leading term in the asymptotic expansion of the
interaction force is (see e.g.~\cite{KK}):

\begin{equation}
\label{force_lub}
F_{1\rightarrow 2}=-\kappa
\frac{u_2-u_1}{d},
\end{equation}
where $\kappa$ is a constant which depends on the viscosity of the lubricating fluid and the radii of the spheres.
We shall take  $\kappa=1$ in what follows.
Consider now (see Fig.~\ref{fig:geom}) an array of $N+1$ spheres, on the $x$-axis, with the same radius $\eps$. We set the first and the last sphere at positions $0$ and $1$, respectively. As a consequence, the number of  degrees of freedom is $N-1$, whereas the number of actual spheres is $N+1$.
\vspace{0.5cm}

\begin{definition}
\label{definition:eps_feasible}
Given a vector of positions $\qq=(q_i)_{1\leq i\leq
N-1}$, we say that $\qq$ is $\eps$-feasible (spheres do not overlap), if 
$$
q_i-q_{i-1}-2\eps \geq 0 \quad \forall i = 1\ldots N, \text{ with }
q_0=0 \text{ and } q_N=1
$$
and strictly $\eps$-feasible if all inequalities are strict (spheres do not touch).
\end{definition}
We denote by  $d_i=q_i-q_{i-1}-2\eps$ the distance between spheres $i$
and $i-1$, by $u_i$ the instantaneous velocity of sphere
$i$, and by $\uu=(u_i)_{1\leq i \leq N-1}$ the velocity vector.
Velocities of the extremal spheres $0$ and $N$ are taken as $0$ (see remark~\ref{remark:nonhom} for non-zero extremal velocities).
Given a strictly $\eps$-feasible vector $\qq$, we define 
$A(\qq)$ as the $(N-1)\times (N-1)$ tridiagonal stiffness matrix
\begin{equation}
\label{A}
 A(\qq)=\left(
\begin{array}{ccccccc}
\frac{1}{d_1}+\frac{1}{d_2} & -\frac{1}{d_2} & & & & & \\
-\frac{1}{d_2} &\frac{1}{d_2}+\frac{1}{d_3} & -\frac{1}{d_3} & & &
& \\
 & \ddots & \ddots & \ddots & & & \\
 & & -\frac{1}{d_i} &\frac{1}{d_i}+\frac{1}{d_{i+1}}
 &-\frac{1}{d_{i+1}}& & \\
 & & &\ddots &\ddots &\ddots & \\
 & & & & -\frac{1}{d_{N-2}} &\frac{1}{d_{N-2}}+\frac{1}{d_{N-1}} &
-\frac{1}{d_{N-2}} \\
 & & & & & -\frac{1}{d_{N-1}} &\frac{1}{d_{N-1}}+\frac{1}{d_N} \\
\end{array}
\right)
\end{equation}


with $d_i=q_i-q_{i-1}-2\eps$.
Consider now a set of forces $f_1$, $f_2$,\dots $f_{N-1}$, and the corresponding vector $\ff$. From~(\ref{force_lub}), the balance of forces  reads

\begin{equation}
\label{eq:sys}
  -A(\qq)\uu+
  \ff = 0.
\end{equation}

\begin{proposition}
\label{proposition:sol_ss_masse}
 Given a strictly $\eps$-feasible vector $\qq\in\IR^{N-1}$, a force
 field $\ff\in\IR^{N-1}$, problem~(\ref{eq:sys}) has a unique solution
 $\uu$, and we shall write
$$
\uu = (u_i)_{1\leq i\leq
N-1}=
 \PP(\qq,\ff,\eps).
$$
This solution can be written
\begin{equation}
\label{eq:sol_sys_disc}
u_i=\frac{1}{D_N}\left\{ (D_N-D_i)\sum_{k=1}^{i}D_k f_k +
D_i\sum_{k=i+1}^{N-1}(D_N-D_k)f_k \right\} \quad \forall i=1\ldots
N-1
\end{equation}
with $\displaystyle D_i=\sum_{j=1}^i d_j$.
\end{proposition}
\begin{proof}
Matrix $A$, which is similar to the matrix obtained by discretizing the Laplace operator with Dirichlet boundary condition by finite differences, is symmetric positive definite, and the vector $\uu$ is immediately checked to solve the system.
\end{proof}

\begin{figure}[th]
\psfragscanon
\psfrag{qim1}[l]{$q_{i_k-1}$}
\psfrag{qi}[l]{$q_{i_k}$}
\psfrag{qiN}[l]{$q_{i_k+N_k}$}
\psfrag{qiNp1}[l]{$q_{i_k+N_k+1}$}
\psfrag{di}[l]{$d_{i_k}$}
\psfrag{diNp1}[l]{$d_{i_k+N_k+1}$}
\psfrag{eps}[l]{$2\varepsilon$}
\centerline{\resizebox{1.0\textwidth}{!}{
{\includegraphics{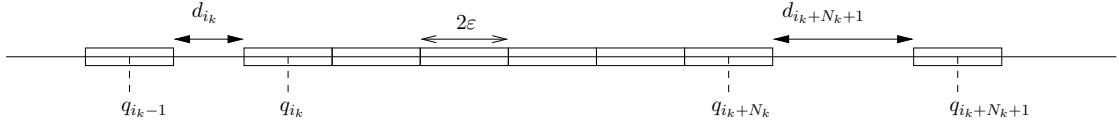}}
}}
\vspace*{8pt}
\caption{Non-strictly $\eps$-feasible configuration}
\label{fig:feas}
\end{figure}

This approach can be extended to $\eps$-feasible situations in a large sense (particles are allowed to get into contact). As the interaction force (which tends to penalize the relative velocity) blows up when particles tend to get into contact, we simply consider that two particles in contact have the same velocity. 
 This situation can be formalized the following way:
 The $N+1$ particules form $N_c$ clusters, and the 
 $k$-th cluster contains the $N_k+1$ particles  $i_k$, 
   $i_k+1$, \dots,  $i_k+N_k$ (see Fig.~\ref{fig:feas}). The balance
 of forces now reads

\begin{eqnarray}
\label{eq:sys_paquets1}
\dsp\forall i \notin \cup_k [i_k,i_k+N_k] \hbox{,      }&
\dsp\frac{u_{i+1}-u_i}{d_{i+1}} - \frac{u_i-u_{i-1}}{d_{i}}=-f_i,\\
\label{eq:sys_paquets2}
\dsp\forall k \in [1,N_c]\hbox{,      }&
\left\{
\begin{array}{l}
\dsp u_{i_k}=u_{i_{k+1}}=\ldots=u_{i_k+N_k}, \\
\dsp\frac{u_{i_k+N_k+1}-u_{i_k+N_k}}{d_{i_k+N_k+1}} -
\frac{u_{i_k}-u_{i_k-1}}{d_{i_k}}=-\sum_{i=i_k}^{i_k+N_k} f_i.
\end{array}
\right.
\end{eqnarray}

\begin{proposition}
\label{proposition:sol_ss_massebis}
 Given an  $\eps$-feasible vector $\qq\in\IR^{N-1}$, a force field $\ff\in\IR^{N-1}$, problem~(\ref{eq:sys_paquets1}-\ref{eq:sys_paquets2}) has a unique solution, and we shall write as before
$$
\uu = (u_i)_{1\leq i\leq
N-1}=
 \PP(\qq,\ff,\eps).
$$
An explicit expression of this solution is given by (\ref{eq:sol_sys_disc}).
\end{proposition}

\begin{remark}
\label{remark:nonhom}
It is possible to set extremal velocities $u_0$ and $u_N$ to non-zero
values: in that case, the balance of forces is given by
\[ A(\qq)\uu = \ff +\bb, \]
where $\bb$ contains the non homogeneous dirichlet boundary conditions:
$$
\bb=(u_0/d_1,0,\ldots,0,u_N/d_N)^t.
$$
 The extension of
Proposition~\ref{proposition:sol_ss_massebis} to that case is straightforward.
\end{remark}

\section{Asymptotic Behaviour of the Discrete Solutions}
\label{sec:micmacstat}



\begin{figure}[th]
\psfragscanon
\psfrag{qim1}[l]{$q_{i_k-1}$}
\psfrag{qi}[l]{$q_{i_k}$}
\psfrag{qiN}[l]{$q_{i_k+N_k}$}
\psfrag{qiNp1}[l]{$q_{i_k+N_k+1}$}
\psfrag{ui}[l]{$u_{i_k}$}
\psfrag{uim1}[l]{$u_{i_k-1}$}
\psfrag{uiNp1}[l]{$u_{i_k+N_k+1}$}
\psfrag{uiN}[l]{$u_{i_k+N_k}=u_{i_k}$}
\psfrag{eps}[l]{$2\varepsilon$}
\centerline{\resizebox{1.0\textwidth}{!}{
{\includegraphics{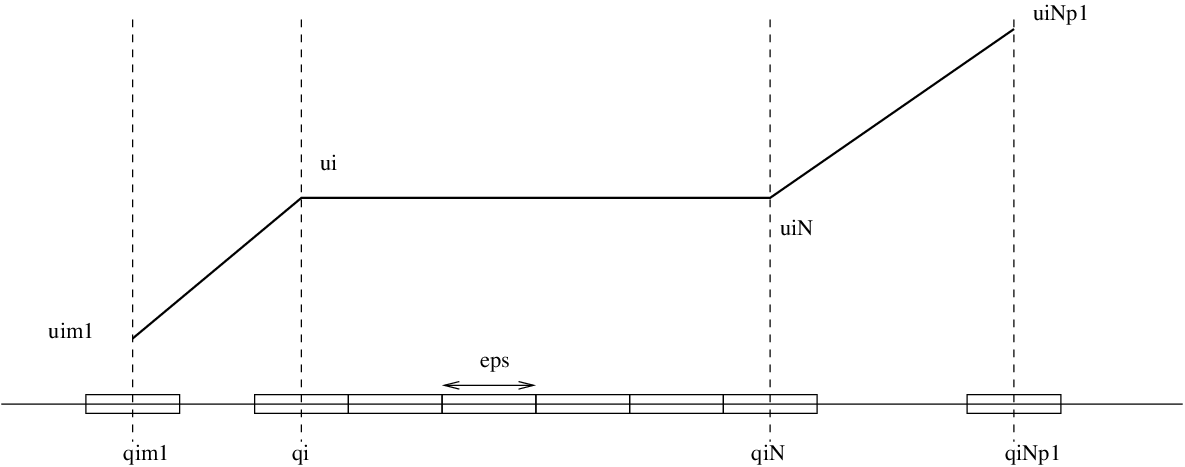}}
}}
\vspace*{8pt}
\caption{Definition of  $u = \tPP(\qq,f,\eps)$}
\label{fig:u}
\end{figure}

Let  $I$ denote   $]0,1[$. 
Firstly, we build  a new operator $\tPP$, which is our key tool to connect the microscopic level to the macroscopic one.
This operator is defined the following way:
 Given $\eps>0$, an $\eps$-feasible position vector $\qq\in\IR^{N-1}$ (as stated by Definition~\ref{definition:eps_feasible}, 
  it represents 
 a distribution of $N+1$ particles whose centers are located at
  $q_0=0$, $q_i$ for $i=1,\dots,N-1$ and $q_N=1$,  with common radius $\eps$),
  a force density $f\in L^1(I)$, we define 
$\uu$ as $\PP(\qq,\ff^\eps,\eps)$ (see
  Proposition~\ref{proposition:sol_ss_masse}
  or~\ref{proposition:sol_ss_massebis}, depending on whether $\qq$ is
  feasible or strictly feasible), where $\ff^\eps$ is defined by
$$
f_i^\varepsilon=\frac{1}{2\varepsilon}\int_{q_i-\varepsilon}^{q_i+\varepsilon}f(s)ds \hbox{ for } 1
\leq i \leq N-1.
$$
Now, to this vector  $\uu=\PP(\qq,\ff^\eps,\eps)=(u_i)_{i=1\ldots N-1}$, we associate a piecewise affine function $u$ defined by
\begin{equation}
\label{eq:udef}
\begin{array}{c}
u \in C^0(\overline I)\virg
 u \hbox{  affine on } [q_i,q_{i+1}] \quad \forall i=0,\dots, N-1\\
u(q_i)=u_i \quad \forall i=1,\dots, N-1 \virg u(0)=u(1)=0.
\end{array}
\end{equation}
We shall write $u = \tPP(\qq,f,\eps)$ (see Fig.~\ref{fig:u}).

In the same spirit, for any $\eps>0$, and any $\eps$-feasible position vector $\qq\in\IR^{N-1}$, we  define $\chi(\qq,\eps)$ as the characteristic function of the solid phase associated to the $(\qq,\eps)$ distribution:
\begin{equation}
\label{eq:chidef}
\chi(\qq,\eps)=\sum_{i=1}^{N-1}\ONE_{]q_i-\varepsilon ,
q_i+\varepsilon[}+\ONE_{]0,
\eps[}+\ONE_{]1-\eps , 1[}
.
\end{equation}

Before we state the main convergence theorem, we still have to give a sense to~(\ref{macro_ss_masse}) when the density $\rho\in [0,1]$ is allowed to take value $1$, even on a set of non-zero measure.

\begin{proposition}
\label{proposition:def_sol_Kinfini}
Let  $K : I \mapsto \IR^+\cup \{+\infty\}$ be measurable, $K(x) \geq \alpha >0$ for almost every $x\in I$,
 $\varphi\in H^{-1}(I)$, and let $J$ be defined as
$$
v \in H^1_0(I) \longmapsto J(v)=\int_IK(x)|\partial_x v|^2 - <\varphi,v>\;\in  \IR\cup \{+\infty\}.
$$
There exists a unique $u\in H^1_0(I)$ which realizes the minimum of $J$ over $H^1_0(I)$.
If there exists  $f\in L^1(I)$  such that $<\varphi,v>=\int fv$, we shall say that $u$ is a generalized solution to
$$
 -\partial_x\left(K(x)\partial_xu\right)=f.
$$
\end{proposition}

\begin{proof}
The functional $J$ is  convex (strictly convex over its domain), coercive, and it can be written
$$
J(v)=\sup_{n\in \IN}\left (\int_I\min(K(x),n)|\partial_x v|^2 - \int_I fv\right ),
$$
thus it is l.s.c. as a supremum of a family of l.s.c
functions. Therefore it admits a unique minimizer.
Note that the minimization problem is equivalent to the problem which consists in minimizing the same functional $J$ over
the set
\begin{equation}
\label{eq:HK}
H_K=\left\{v\in H^1_0(I), \partial_xv=0 \hbox{ a.e. in } D(K)^c
\hbox{ and } \int_{D(K)}K(x)|\partial_x v|^2<\infty\right\}
.
\end{equation}
where $D(K)=\{x\in I, K(x)<+\infty\}$ is the domain of $K$.
 Consequently, $u$ is  characterized by the variational formulation
\begin{equation}
\label{eq:form_var}
u\in H_K \virg \int_{D(K)}K(x)\partial_x u\,\partial_x v = \int_I
fv\quad \forall v\in H_K.
\end{equation}
\end{proof}

We may now state the convergence result.
\begin{theorem}
\label{thm:stat}
Let $f\in L^\infty(I)$ be a Lipschitz force density, and $\rho\in L^\infty(I)$ a
solid fraction, with $\rho(x) \in [0,1]$ {a.e.} in $I$. 
Let $(\qq^\eps)_\eps$ be a sequence of $\eps$-feasible position
vectors (see Definition~\ref{definition:eps_feasible}),  
 $ \qq^\eps \in
\IR^{N^\eps-1}$ with $N^\eps=1/\eps$.
We introduce $\chi^\eps = \chi (\qq^\eps,\eps)$ (see~(\ref{eq:chidef})),
and we assume that $\chi^\eps$ converges toward $\rho$ in $L^\infty(I)$ weak-$\star$.

Then $u^\eps=\tPP (\qq^\eps,f,\eps)$ solution to the discrete model (see~(\ref{eq:udef})) converges weakly toward $u$ in $H^1_0(I)$ as
    $\eps$ converges to $0$, where $u$ is the solution to 
\begin{equation}
\label{eq_force}
-\partial_x\left(\frac{1}{1-\rho}\partial_xu\right)=\rho f,
\end{equation}
in the sense of Proposition~\ref{proposition:def_sol_Kinfini}
(i.e. characterized by~(\ref{eq:form_var})).

\end{theorem}

\begin{remark}
\label{rmk:expl_construction}
For any $\rho\in L^\infty(I)$ with $\rho(x) \in [0,1]$ {a.e.} in $I$, one can easily construct a sequence $(\qq^\eps)_\eps$,
$\eps$-feasible, such that  $\chi^\eps\wsto\rho$ in
$L^\infty(I)$. Let $\eps$ be equal to $1/N$. We define
$q_0^\eps=0$. Then, we denote by $l_1^\eps$ the real such that
$\int_{q_0^\eps}^{l_1^\eps}\rho=\eps$. Since $\rho\leq 1$ we have
$l_1^\eps\geq q_0^\eps+\eps$ and we can define
$q_1^\eps=l_1^\eps+\eps$. Similarly, to define $q_1^\eps$, we first
chose $l_2^\eps$ such that
$\int_{l_0^\eps}^{l_1^\eps}\rho=2\eps$ and define
$q_2^\eps=l_2^\eps+\eps$. We iterate to obtain  $(\qq^\eps)_\eps$,
$\eps$-feasible. Moreover we have
$\int_{l_i^\eps}^{l_{i+1}^\eps}(\chi^\eps-\rho)=0$ and it follows that $\chi^\eps\wsto\rho$ in
$L^\infty(I)$ when $\eps$ goes to zero (see the proof of
lemma~\ref{rho_cv}, where a similar property is established).
\end{remark}

\proof
The proof is based on some technical lemmas. 
For readability reasons, we postpone the proofs of the lemmas to the end of the section.

\begin{figure}[th]
\psfragscanon
\psfrag{qim1}[l]{$q_{i_k-1}^\varepsilon$}
\psfrag{qi}[l]{$q_{i_k}^\varepsilon$}
\psfrag{qiN}[l]{$q_{i_k+N_k}^\varepsilon$}
\psfrag{qiNp1}[l]{$q_{i_k+N_k+1}^\varepsilon$}
\psfrag{rhoi}[l]{$\rho_{i_k}^\varepsilon$}
\psfrag{rho1}[l]{$\rho^\varepsilon=1$}
\psfrag{rhoiNp1}[l]{$\rho_{i_k+N_k+1}^\varepsilon$}
\centerline{\resizebox{0.95\textwidth}{!}{
{\includegraphics{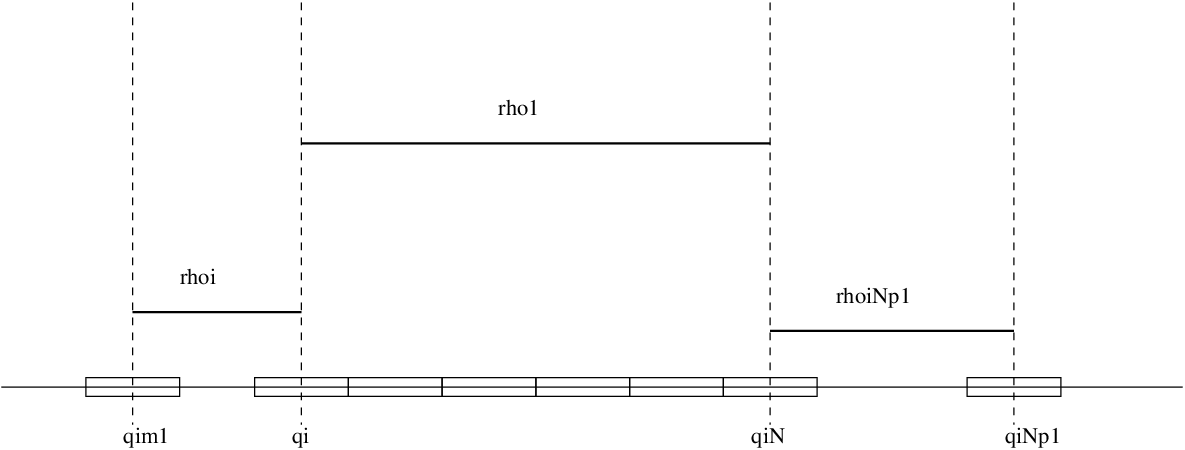}}
}}
\vspace*{8pt}
\caption{Definition of $\rho^\eps$}
\label{fig:rho_eps}
\end{figure}

As a first step, we define  $\rho^\eps$, piecewise constant, as the proportion of solid on each
subinterval $[q_{i-1}^\eps,q_i^\eps]$ in the following way: let $d_i^\eps$ be the distance between particles $i-1$ and $i$, then
\begin{equation}
\label{eq_rho_eps}
\forall i=1,\ldots,N^\eps \virg \rho^\varepsilon=\rho^\varepsilon_i=1-\frac{d_{i}^\varepsilon}{q_{i}^\varepsilon-q_{i-1}^\varepsilon}
\hbox{ on } [q_{i-1}^\varepsilon,q_{i}^\varepsilon].
\end{equation}
where $q_0^\eps=0$ and $q_{N^\eps}^\eps=1$. This definition is
illustrated by Fig.~\ref{fig:rho_eps}. Note that, in this figure and in the
following of the proof, we shall drop the superscript $\eps$ for
indexes concerning the clusters in order to alleviate notations (so that,
for instance, we shall write $N_c$ for $N_c^\eps$,
$q_{i_k}^\varepsilon$ for $q_{i_k^\eps}^\varepsilon$ or
$q_{i_k+N_k}^\varepsilon$ for $q_{i_k^\eps+N_k^\eps}^\varepsilon$).

\begin{figure}[th]
\psfragscanon
\psfrag{a}[l]{$a=q_{0}^\varepsilon$}
\psfrag{q1}[l]{$q_{1}^\varepsilon$}
\psfrag{qim1}[l]{$q_{i_k-1}^\varepsilon$}
\psfrag{qi}[l]{$q_{i_k}^\varepsilon$}
\psfrag{qiN}[l]{$q_{i_k+N_k}^\varepsilon$}
\psfrag{qiNp1}[l]{$q_{i_k+N_k+1}^\varepsilon$}
\psfrag{qNm1}[l]{$q_{N^\eps-1}^\varepsilon$}
\psfrag{b}[l]{$b=q_{N^\eps}^\varepsilon$}
\psfrag{w1}[l]{$w_{1}^\varepsilon$}
\psfrag{w2}[l]{$w_{2}^\varepsilon$}
\psfrag{wim1}[l]{$w_{i_k-1}^\varepsilon$}
\psfrag{wi}[l]{$w_{i_k}^\varepsilon$}
\psfrag{wiNp1}[l]{$w_{i_k+N_k+1}^\varepsilon$}
\psfrag{Bip1}[l]{$w_{i_k+1}^\varepsilon=\beta_{i_k+1}^\varepsilon$}
\psfrag{BiN}[l]{$w_{i_k+N_k}^\varepsilon=\beta_{i_k+N_k}^\varepsilon$}
\psfrag{wNm1}[l]{$w_{N^\eps-1}^\varepsilon$}
\psfrag{wN}[l]{$w_{N^\eps}^\varepsilon$}
\centerline{\resizebox{\textwidth}{!}{
{\includegraphics{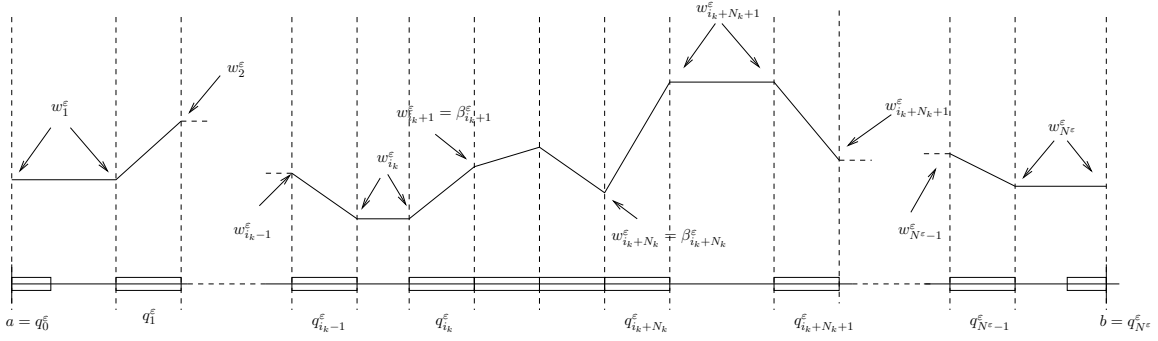}}
}}
\vspace*{8pt}
\caption{Definition of $w^\varepsilon$}
\label{fig:w_eps}
\end{figure}

%
\flushleft Define now $w^\varepsilon$ (see Fig.~\ref{fig:w_eps}),
affine by part, as the discrete counterpart of $\partial_xu/(1-\rho)$ given by
\begin{equation}
\label{eq_w_eps}
\left\{
\begin{array}{ll}
w^\varepsilon(q_i^\varepsilon-\varepsilon)=w_i^\varepsilon&\hbox{ for }
i=1\ldots N^\varepsilon -1\vseq\\
w^\varepsilon(q_i^\varepsilon+\varepsilon)=w_{i+1}^\varepsilon
\end{array}
\right.
\end{equation}
with
\begin{equation}
\label{eq_wi_eps}
\left\{
\begin{array}{ll}
w^\varepsilon_i=\dsp\frac{1}{1-\rho^\varepsilon_i}\partial_xu^\varepsilon&
\hbox{, if  } \rho_i^\varepsilon < 1 \\
w^\varepsilon_i=\beta_i^\eps & \hbox{, if  } \rho_i^\varepsilon =  1\\
\end{array}
\right.
\end{equation}
where $(\beta_i^\eps)_{i_k+1\leq i \leq i_k+N_k}$ corresponding to the
  kth cluster, is the solution to the following system:
\begin{equation}
\label{eq:sys_beta}
\left\{
\begin{array}{ccccc}
\beta_{i_k+1}^\varepsilon&-&\dsp\frac{u_{i_k}^\varepsilon-u_{i_k-1}^\varepsilon}{d_{i_k}^\varepsilon}&=&-2\varepsilon f_{i_k}^\varepsilon,\vspace{6pt}\\
\beta_{i_k+2}^\varepsilon&-&\beta_{i_k+1}^\varepsilon&=&-2\varepsilon f_{i_k+1}^\varepsilon,\\
 & & & \vdots & \\
\beta_{i_k+N_k}^\varepsilon&-&\beta_{i_k+N_k-1}^\varepsilon&=&-2\varepsilon
f_{i_k+N_k-1}^\varepsilon,\vspace{6pt}\\
\dsp\frac{u_{i_k+N_k+1}^\varepsilon-u_{i_k+N_k}^\varepsilon}{d_{i_k+N_k+1}^\varepsilon}&-&\beta_{i_k+N_k}^\varepsilon&=&-2\varepsilon
f_{i_k+N_k}^\varepsilon.\\
\end{array}
\right.
\end{equation}

Note that, summing up all equations of~(\ref{eq:sys_beta}) we recognize the
balance of forces on the kth cluster given by (\ref{eq:sys_paquets2}).

\begin{remark}
The idea behind the above construction is that $\beta_i$ can be seen as the cohesion
  force between particles $i-1$ and $i$. 
A first way to notice it is to note that $\uueps$ is the limit of
  $\uu^{\eps,\eta}$  where $\uu^{\eps,\eta}$ is the solution to
  system~(\ref{eq:sys}) with $d_i=\eta>0$ for $i$ between $i_k+1$ and
  $i_k+N_k$ and that we have 
\[\dsp\forall k \virg \forall j \in [1,N_k] \virg \beta_{i_k+j}^\eps=\lim_{\eta\rightarrow 0}
  \frac{u_{i_k+j}^{\eps,\eta}-u_{i_k+j-1}^{\eps,\eta}}{\eta} .\] 
Another way to define these cohesion forces is to consider the
following 
minimization problem: minimize
\[\dsp J(\vv)=\sum_{i\notin \cup_k [i_k+1,i_k+N_k]}
    \frac{1}{2}\frac{(v_i-v_{i-1})^2}{d_i^\eps} + \sum_{i=1}^{N-1}
    f_i^\eps v_i\]
over $\dsp K=\{\vv \virg \forall i \in \cup_k [i_k+1,i_k+N_k] ,\,
  v_i=v_{i-1}\}$. This problem is equivalent to
  problem~(\ref{eq:sys_paquets1}-\ref{eq:sys_paquets2}) and
$\beta_i^\eps$ turns up as the Lagrange multiplier
associated to the constraint $v_i=v_{i-1}$.
\end{remark}

The next lemma shows that~(\ref{eq_force}) is true at the discrete level:

\begin{lemma}
\label{lem_u_eps} For any $\eps>0$,
\begin{equation}
\label{eq_force_eps}
-\partial_x\left(w^\varepsilon\right)=
f^\varepsilon
\end{equation}
where $f^\varepsilon={\dsp
\sum_{i=1}^{N^\varepsilon-1}}f_i^\varepsilon {\mathbf 1}_{]q_i^\varepsilon-\varepsilon,q_i^\varepsilon+\varepsilon[}.$
\end{lemma}

The idea of the proof is now to let $\eps$ go to $0$ in~(\ref{eq_force_eps}). To that purpose, we first study
$\rho^\varepsilon$ and $2\varepsilon f^\varepsilon$ when $\eps$ tends
to zero:

\begin{lemma}
\label{rho_cv} $\rho^\varepsilon \wsto \rho$ in $L^\infty(I)$.
\end{lemma}

\begin{lemma}
\label{f_cv}
 $f^\varepsilon \wsto \rho f$ in $L^\infty(I)$.
\end{lemma}


To study the convergence of $(u^\eps)_\eps$, we prove in next lemma
that it is bounded in $H^1_0(I)$.

\begin{lemma}
\label{u_eps_bornee}
 $(u^\eps)_\eps$ is bounded in $H^1_0(I)$.
\end{lemma}

Hence, we can extract a subsequence of $(u^\eps)_\eps$
    (still denoted by $(u^\eps)_\eps$ ) such that
$$ u^\eps \rightharpoonup u \hbox{ in
    } H^1_0(I).$$
In order to pass to the limit in the left-hand side
of~(\ref{eq_force_eps}), we are going to study the convergence of
$(w^\eps)_\eps$. It will follow from the next lemma.
\begin{lemma}
\label{w_eps_bornee}
$(w^\varepsilon)_\varepsilon$ and $(\partial_x w^\varepsilon)_\varepsilon$ are bounded in $L^\infty(I)$.
\end{lemma}
Consequently, from Ascoli theorem, we can find $w$ in ${\mathcal C}^0(I)$ and a subsequence of
$(w^\eps)_\eps$ (still denoted by $(w^\eps)_\eps$ ) such that
\begin{equation} 
\label{w_eps_cv} 
w^\varepsilon\rightarrow w \hbox{ in } {\mathcal C}^0(I).
\end{equation}
Then, convergence of $w^\eps$ and $\rho^\eps$ make it possible to establish the following lemma:
\begin{lemma}
\label{dx_u_eps}
$\partial_x u^\varepsilon \wsto (1-\rho)w \hbox{ in
} L^\infty(I).$
\end{lemma}
We now come to the last step of the proof: 
let $\varepsilon$ tend
to zero in~(\ref{eq_force_eps}) and obtain asymptotically
\[
-\partial_x\left(\frac{1}{1-\rho}\partial_xu\right)=\rho
f,
\]
in the sense of Proposition~\ref{proposition:def_sol_Kinfini}, characterized by~(\ref{eq:form_var}).

First, Lemma~\ref{dx_u_eps} implies $\partial_xu=(1-\rho)w$, so 
$\dsp u\in H_{1/(1-\rho)}$ (defined by~(\ref{eq:HK})).
Then,  by Lemma~\ref{lem_u_eps},
\[ \int_Iw^\varepsilon v' = \left< 2\varepsilon f^\varepsilon ,
v \right> \hbox{   } \forall v\in H^1_0(I).\]
Passing to the limit on $\eps$ gives, using~(\ref{w_eps_cv}) and Lemma~\ref{f_cv},
\[ \int_I w v' = \int_I\rho f v \hbox{   }\quad \forall v\in H^1_0(I).
\]
Finally, for any $\dsp v\in H_{1/(1-\rho)}$, 
\[
\begin{array}{l}
\displaystyle\int_I \rho fv = \int_I wv'= \int_{\rho \neq 1}wv'\vspace{6pt}\\
\displaystyle\hspace{1cm} = \int_{\rho \neq 1} \frac{(1-\rho)w}{1-\rho}v' =
\int_{\rho \neq 1} \frac{\partial_x u}{1-\rho}v'
\end{array}
\]
and we conclude that $u$ is the solution to (\ref{eq_force}).

So, we proved that there exists a subsequence of $(u^\eps)_\eps$ converging 
to $u$ as $\eps$ tends to zero. Since the same work can be done for
each subsequence of $(u^\eps)_\eps$, we conclude that
$(u^\eps)_\eps$ itself converges to $u$, which completes the proof of
the theorem.
\qed

\bigskip

{\bf 
{Proof of Lemma \ref{lem_u_eps}}}: $-\partial_x\left(w^\varepsilon\right)=2\varepsilon
f^\varepsilon$.

First, a simple computation shows that
\[
\begin{array}{l}
\dsp \partial_x w^\varepsilon =
\sum_{\tiny \begin{array}{c}\tiny i\in
  1\ldots N^\varepsilon-1\\\tiny d_i>0\hbox{,  }d_{i+1}>0\end{array}}\frac{1}{2\eps}\left(\frac{\ueps_{i+1}-\ueps_i}{d_{i+1}^\varepsilon}-\frac{\ueps_{i}-\ueps_{i-1}}{d_{i}^\varepsilon}\right)\ONE_{]q_i^\varepsilon-\varepsilon,q_i^\varepsilon+\varepsilon[}
\vspace{6pt}\\
\dsp\quad\quad
+\sum_{k=1}^{P} \sum_{i=i_k+1}^{i_k+N_k-1} \frac{1}{2\eps}\left(\beta_{i+1}^\eps-\beta^\eps_{i}\right)\ONE_{]q_i^\varepsilon-\varepsilon,q_i^\varepsilon+\varepsilon[}
\vspace{6pt}\\
\dsp\quad\quad
+\sum_{k=1}^{P} \frac{1}{2\eps}\left(
\beta^\eps_{i_k+1}-\frac{\ueps_{i_k}-\ueps_{i_{k-1}}}{d_{i_k}^\eps}
\right)\ONE_{]q_{i_k}^\varepsilon-\varepsilon,q_{i_k}^\varepsilon+\varepsilon[}
\vspace{6pt}\\
\dsp\quad\quad
+\sum_{k=1}^{P} \frac{1}{2\eps}\left(
\frac{\ueps_{i_k+N_k+1}-\ueps_{i_k+N_k}}{d_{i_k+N_k+1}^\eps}-\beta^\eps_{i_k+N_k}
\right)\ONE_{]q_{i_k+N_k}^\varepsilon-\varepsilon,q_{i_k+N_k}^\varepsilon+\varepsilon[}.\vspace{6pt}
\end{array}
\]
Then, combining this with system
(\ref{eq:sys_paquets1}-\ref{eq:sys_paquets2}) and with the
definition of $\beta_i$~(\ref{eq:sys_beta}), we get
$$ \partial_x w^\varepsilon
=-\sum_{i=1}^{N^\varepsilon-1}
f_i^\varepsilon\mathbf{1}_{]q_i^\eps-\eps,q_i^\eps+\eps[}.
$$
\qed

{\bf
{Proof of Lemma \ref{rho_cv}}}: $\rho^\varepsilon \wsto \rho$ in $L^\infty(I)$.

Since $\rho^\eps-\rho=(\rho^\eps-\chi^\eps)+(\chi^\eps-\rho)$ and
$\chi^\eps \wsto \rho$ in $L^\infty(I)$, the result will follow
    provided we show that
\[\forall \varphi \in  L^1(I) \virg
\lim_{\varepsilon\rightarrow 0} \left( \int_I \chi^\varepsilon
\varphi - \int_I \rho^\eps\varphi \right) =0.\]

By density of the set of stairs functions in $L^1(I)$, and
    by using the fact that $(\rho^\varepsilon)_\varepsilon$ and
    $(\chi^\varepsilon)_\varepsilon$ are bounded in $L^\infty(I)$,
    this in turn will follow from
   
$$
\forall \varphi \hbox{ piecewise constant on  } I \virg
\lim_{\varepsilon\rightarrow 0} \left( \int_I \chi^\varepsilon
\varphi - \int_I \rho^\eps\varphi \right) =0.$$
To show this, it suffices to prove that
\[\forall \alpha,\beta \virg 0<\alpha<\beta<1 \virg \lim_{\varepsilon\rightarrow
0}\left(
\int_\alpha^\beta\chi^\varepsilon-\int_\alpha^\beta\rho^\varepsilon
\right) =0.
\]
In order to do so take $\alpha$ and $\beta$ such that
$0<\alpha<\beta<1$ and denote the particules whose centers are in
$[\alpha,\beta]$ by $q_{i_0}^\varepsilon$, $q_{i_0+1}^\varepsilon$, \ldots,
$q_{j_0}^\varepsilon$ (again we drop the $\eps$ superscript,
keeping in mind that $i_0$ and $j_0$ depend on $\eps$). Since $\int_{q_i^\varepsilon}^{q_{i+1}^\varepsilon}(\chi^\varepsilon-\rho^\varepsilon)=2\varepsilon-(q_{i+1}^\varepsilon-q_i^\varepsilon-d_{i+1}^\varepsilon)=0$
for $1\leq i\leq N^\varepsilon-1$, we have 
\[\int_\alpha^\beta\chi^\varepsilon-\int_\alpha^\beta\rho^\varepsilon=\int_\alpha^{q_{i_0}^\varepsilon}(\chi^\varepsilon-\rho^\varepsilon)-\int_{q_{j_0}^\varepsilon}^\beta(\chi^\varepsilon-\rho^\varepsilon).
\]
Then, a simple computation shows that the left-hand side converges to zero
as $\eps$ tends to zero wich completes the proof of the lemma.
\qed

{\bf 
{Proof of Lemma~\ref{f_cv}}}: $f^\varepsilon \wsto \rho f$ in
$L^\infty(I)$.

%
Writing  
\[f^\eps - \rho f = \left( f^\eps - \chi^\eps f\right) + \left( \chi^\eps f - \rho f\right).
\]
and using the fact that $\chi^\varepsilon \wsto \rho$ in $L^\infty(I)$, the required result will follow as soon as we prove that 
\[ \forall \varphi \in L^1(I) \virg
\lim_{\eps\rightarrow 0} \int_I  \left( f^\eps - \chi^\eps f\right)\varphi = 0,\]
To obtain this, merely compute
\[
\begin{array}{rcl}
\displaystyle \int_I  \left( f^\eps - \chi^\eps f\right)\varphi &=&\displaystyle \sum_{i=1}^{N^\varepsilon-1}
\int_{q_i^\varepsilon-\varepsilon}^{q_i^\varepsilon+\varepsilon}[f_i^\eps-f(x)]\varphi(x)dx
\vspace{6pt}\\
&&\displaystyle-\int_{0}^{\varepsilon}f(x)\varphi(x)dx
-\int_{1-\varepsilon}^{1}f(x)\varphi(x)dx
\end{array}
\]
and, since the last two terms tend to zero as $\eps$ go to zero, the result will follow as soon as we prove
\[
\lim_{\eps\rightarrow 0} \sum_{i=1}^{N^\varepsilon-1}
\int_{q_i^\varepsilon-\varepsilon}^{q_i^\varepsilon+\varepsilon}[f_i^\eps-f(x)]\varphi(x)dx =0.
\]
To do so, use the Lipschitz hypothesis on $f$ to write
\[
\left|f_i^\eps-f(x)\right|=
\frac{1}{2\eps}\left|\int_{q_i^\varepsilon(t)-\varepsilon}^{q_i^\varepsilon(t)+\varepsilon}
(f(y,t)-f(x,t))dy\right|\leq 2C\eps,
\]
which gives
\[
\left|\int_0^1 (f^\eps-\chi^\eps f)\varphi \right|\leq
2C\eps \|\varphi\|_{L^1(I)}
\]
and completes the proof of the lemma.
\qed

{\bf 
{Proof of Lemma~\ref{u_eps_bornee}}}: $(u^\varepsilon)_\varepsilon$  is
bounded in $H^{1}_0(I)$.

Using the fact that $\partial_x u^\eps$ is zero when
$\rho^\eps$ is equal to $1$ we obtain
\[
\|u^\varepsilon\|_{H^1_0(I)}^2=
\int_{I\cap
  \{\rho^\eps<1\}}|\partial_xu^\varepsilon|^2\leq
\int_{I\cap
  \{\rho^\eps<1\}}\frac{1}{1-\rho^\eps}|\partial_x
u^\eps|^2
\]
and
\begin{equation}
\label{eq:dxueps}
\|u^\varepsilon\|_{H^1_0(I)}^2\leq\int_{I\cap \{\rho^\eps<1\}}\tilde
w^\eps\partial_xu^\varepsilon=
\int_{I}\tilde
w^\eps\partial_xu^\varepsilon,
\end{equation}
where $\tilde w^\eps$ is the function, piecewise constant, with $\tilde
w^\eps=w_i^\eps$ on $[q_{i-1}^\varepsilon(t),q_{i}^\varepsilon(t)]$.

A simple computation shows that, in the sense of distribution, we have
\[\displaystyle \partial_x \tilde w^\varepsilon =
\sum_{\scriptsize\begin{array}{c}i\in
  1..N^\varepsilon-1\\d_i^\eps>0\hbox{,  }d_{i+1}^\eps>0\end{array}}\left(\frac{\ueps_{i+1}-\ueps_i}{d_{i+1}^\varepsilon}-\frac{\ueps_{i}-\ueps_{i-1}}{d_{i}^\varepsilon}\right)\delta_{\qeps_i}
+\sum_{k=1}^{N_c} \sum_{i=i_k+1}^{i_k+N_k-1} \left(\beta_{i+1}^\eps-\beta^\eps_{i}\right)\delta_{\qeps_i}
\]
\[\dsp
+\sum_{k=1}^{N_c} \left(
\beta^\eps_{i_k+1}-\frac{\ueps_{i_k}-\ueps_{i_{k-1}}}{d_{i_k}}
\right)\delta_{\qeps_{i_k}}
+\sum_{k=1}^{N_c} \left(
\frac{\ueps_{i_k+N_k+1}-\ueps_{i_k+N_k}}{d_{i_k+N_k+1}}-\beta^\eps_{i_k+N_k}
\right)\delta_{\qeps_{i_k+N_k}}
\]
Then, combining this with system
(\ref{eq:sys_paquets1}-\ref{eq:sys_paquets2}) and with the
definition of $\beta_i$~(\ref{eq:sys_beta}), we get
$$ \partial_x \tilde w^\varepsilon
=-\sum_{i=1}^{N^\varepsilon-1}2\varepsilon
f_i^\varepsilon\delta_{q_i^\varepsilon},
$$
which can be extended to $H^{-1}(I)$ by density. This,
with~(\ref{eq:dxueps}) and the fact that $\qqeps$ is $\eps$-feasible, gives
\[
\displaystyle\|u^\varepsilon\|_{H^1_0(I)}^2\leq
\sum_{i=1}^{N^\eps-1}2\varepsilon 
f_i^\varepsilon u^\eps(q_i^\varepsilon)\leq 
\|u^\varepsilon\|_{L^\infty(I)}\|f\|_{L^1(I)}.
\]
and, by continuous injection of $H^1_0(I)$ in ${\mathcal C}^0(I)$ we
obtain
\[
\|u^\varepsilon\|_{H^1_0(I)}\leq
C\|f\|_{L^1(I)}
\]
which completes the proof of the lemma.
\qed

{\bf 
{Proof of Lemma~\ref{w_eps_bornee}}}: $(w^\varepsilon)_\varepsilon$ and
$(\partial_x w^\varepsilon)_\varepsilon$ are bounded in $L^\infty(I)$.

By~(\ref{eq_w_eps}-\ref{eq_wi_eps}), it suffices to obtain an upper-bound for $(w_i^\eps)_i$, where
$$
w_i^\eps=\frac{1}{1-\rho^\varepsilon_i}\partial_xu^\varepsilon=\frac{u_i^\eps-u_{i-1}^\eps}{d_i^\eps}\hbox{ if }  \rho_i^\varepsilon < 1,
$$
 and $w_i^\eps=\beta_i^\eps$ otherwise.
Note that, by (\ref{eq:sys_beta}) and the fact that $\qqeps$ is $\eps$-feasible
\[\forall k=1\ldots N_c \virg \forall j=1\ldots N_k \virg \left|\beta^\eps_{i_k+j}\right|=\left|w^\eps_{i_k}-2\eps\sum_{m=i_k}^{i_k+j-1}f_m\right|
\leq \left|w^\eps_{i_k}\right|+\|f\|_{L^1(I)}.
\]
Therefore, to show that $\|w^\varepsilon\|_{L^\infty(I)}$ is bounded
it suffices to prove that $(w_i^\varepsilon)_{\{\tiny i\hbox{ s.t. }\rho_i<1\}}$
is bounded. In order to do so, a simple computation gives using~(\ref{eq:sol_sys_disc})
\begin{equation}
\label{wi_eps}
 \forall i\hbox{ s.t. }\rho_i<1 \virg w_i^\varepsilon=2\varepsilon{\dsp
\sum_{k=i+1}^{N^\varepsilon-1}}\frac{D_{N^\varepsilon}-D_k}{D_{N^\varepsilon}}f_k^\varepsilon-2\varepsilon{\dsp
\sum_{k=1}^i}\frac{D_k}{D_{N^\varepsilon}}f_k^\varepsilon.
\end{equation}
Then, from $\frac{D_{N^\varepsilon}-D_k}{D_{N^\varepsilon}}\leq 1$,
$\frac{D_k}{D_{N^\varepsilon}}\leq 1$ and the fact that $\qqeps$ is $\eps$-feasible it follows that 
\[
|w_i^\eps|\leq
\sum_{k=1}^{N^\varepsilon-1}|2\varepsilon f_k^\varepsilon| \leq
\|f\|_{L^1(I)}
\]
and we conclude that
\[ \|w^\varepsilon\|_{L^\infty(I)}\leq
2\|f\|_{L^1(I)}.
\]
The bound on $(\partial_x w^\eps)_\eps$ is easily obtained using
lemma~\ref{lem_u_eps} combined with $f\in L^\infty(I)$, which completes the proof of the lemma.
\qed


%

{\bf 
{Proof of Lemma \ref{dx_u_eps}}}: $\partial_x u^\varepsilon \wsto (1-\rho)w$ in $L^\infty(I)$.

Writing 
\[ \partial_x u^\varepsilon - (1-\rho)w = 
\left\{\partial_x u^\eps - (1-\rho^\eps)w^\eps\right\} + 
\left\{(1-\rho^\varepsilon)(w^\varepsilon-w) \right\} +
\left\{\left((1-\rho^\varepsilon)-(1-\rho)\right)w\right\},
\]
we shall prove the weak-star convergence to zero in $L^\infty(I)$ of
each term of the right-hand side.

The first term goes strongly to zero in $L^\infty(I)$. Indeed, using
the definitions of $w^\eps$~(\ref{eq_w_eps}), $\rho^\eps$~(\ref{eq_rho_eps})
and $u^\eps$~(\ref{eq:udef}) we obtain
\[
\|\partial_x u^\varepsilon -
(1-\rho^\eps) w^\eps\|_\Linf\leq\frac{1}{2}\sup_i|w^\eps_{i+1}-w^\eps_i|.
\]
This, combined with $|w^\eps_{i+1}-w^\eps_i| = |2\eps f^\eps_i|$ gives
\[
\|\partial_x u^\varepsilon -
(1-\rho^\eps)w^\eps\|_\Linf\leq\eps\|f\|_\Linf,
\]
which proves the uniform convergence to zero.

To prove the weak-star convergence of the second term, by density of
  $\mathcal{C}_0^\infty (I)$ in $L^1(I)$ and lemma~\ref{w_eps_bornee},
  it suffices to take test-functions $\varphi$ in
  $\mathcal{C}_0^\infty (I)$. Therefore, the result
    follows immediately from

\[ \left| \int_I (1-\rho^\varepsilon)(w^\varepsilon-w)\varphi \right|
\leq |I| \|\varphi\|_{L^\infty(I)}\|w^\varepsilon-w\|_{L^\infty(I)}\]
together with (\ref{w_eps_cv}).

We shall now prove the weak-star convergence of the last term to
    zero. To do so, take $\varphi$ in $L^1(I)$. Using that $w\in
    {\mathcal C}^0(I)$, it follows that $w\varphi\in
    L^1(I)$ and by Lemma \ref{rho_cv}

\[ \lim_{\varepsilon\rightarrow 0} \int_I \left(
(1-\rho^\varepsilon)-(1-\rho)\right)w\varphi = 0,\]
as required.\qed


\section{Non-stationary model}
\label{sec:micmacnonstat}

In this section, we consider the non-stationnary lubrication model
without inertia. We suppose that the extremal particles are fixed
($q_0(t)=0$, $q_N(t)=1$, $\forall t$). In that case, the unknowns are
the positions of the
particles  $1$ to $N-1$ against time. We denote by
$\qq(t)=(q_1(t),\ldots,q_{N-1}(t))$ these positions. Since the balance of forces is
achieved at each instant $t$, $\qq$ is solution to
\begin{equation}
\label{EDO}
\left|
\begin{array}{l}
\dsp -A(\qq)\dot \qq+ \ff(t,\qq)=0,\vspace{6pt}\\
\dsp \qq(0)=\qqinit,
\end{array}
\right.
\end{equation}
where $A$ is defined in~(\ref{A}) and $\ff=(f_i)_i$ is the vector made
of the external forces
exerted on the particles.


\begin{proposition}
Suppose that $f_i \in L^1_{loc}(\IR_+,L^\infty(\IR^N))$ and is
Lipschitz with respect to the second variable. Let
$\qqinit$ be a strictly $\eps$-feasible vector (see
definition~\ref{definition:eps_feasible}). Then, the system of
ODE~(\ref{EDO}) has a unique global solution.
\end{proposition}
\begin{proof}
Equation~(\ref{EDO}) can be written $\dot\qq=F(t,\qq)$ with
$F(t,\qq)=A^{-1}(\qq)\ff(t,\qq)$. We apply Cauchy-Lipschitz theorem to
this first order ODE on the set of strictly $\eps$-feasible
configurations
$$
\Omega_\eps=\left\{ \qq \virg  d_i = q_i-q_{i-1}-2\varepsilon >0 \hbox{ } \forall
i\in\{1 \ldots N\} \right\}.
$$
If the maximal solution were defined over $[0,T^*[$ with $T^*<+\infty$, 
%
we could find an index $i$ and a
subsequence $(t_n)_n$ such that 
$$
 d_i(t_n) \rightarrow 0 \hbox{ when } n\rightarrow
+\infty \hbox{ and }t_n\rightarrow T^*.
$$
The extremal spheres being fixed, there also exists $k$ such that $d_k(t_n)
\geq \eta>0$ for all $n$. By summing up lines $i$ to $k-1$ of
equation~(\ref{EDO}) and integrating it over $[0,t_n]$ we
obtain
$$
C + \sum_{j=i}^{k-1} \left (
\ln(d_k(t_n))-\ln(d_i(t_n)) \right )= - \sum_{j=i}^{k-1} \int_0^{t_n} f_j(s,\qq(s))ds,
$$
where $C$ is a constant. This is a contradiction since the left-hand
side goes to infinity when $n\rightarrow +\infty$ while the right-hand
side is bounded. Consequently, $T^*=+\infty$ and the solution is global.
\end{proof}

Similarly to what has been done in the previous section, we can
extend this approach to $\eps$-feasible situations in the large
sense. We denote by $\uu$ the velocity $\dot\qq$ of the particles.
\begin{proposition}
\label{prop:Q}
Suppose that $f_i \in L^1_{loc}(\IR_+,L^\infty(\IR^N))$ and is
Lipschitz with respect to the second variable. Let
$\qqinit$ be an $\eps$-feasible vector. The system of
ODE~(\ref{eq:sys_paquets1}-\ref{eq:sys_paquets2}) with initial
condition $\qq(0)=\qqinit$ has a unique global
solution $t\rightarrow \qq(t)$. This solution will be denoted by
$\qq=\QQ(\qqinit,\ff,\eps)$.  The velocity $\uu=\dot\qq$ can be written with respect
to $\qq$ using~(\ref{eq:sol_sys_disc}).
\end{proposition}

The continuous counterpart of this non-stationnary model is the
following system of PDEs made of the constitutive law derived in the
previous section and the transport equation:
$$
\left|
\begin{array}{l}
\displaystyle -\partial_x\left(\frac{1}{1-\rho}\partial_xu\right)=\rho f\vspace{6pt}\\
\displaystyle \partial_t \rho+\partial_x(\rho u)=0
\end{array}
\right.
$$

\begin{figure}[hbtp]
\psfragscanon
\psfrag{qim1}[l]{$q_{i_k-1}(t)$}
\psfrag{qi}[l]{$q_{i_k}(t)$}
\psfrag{qiN}[l]{$q_{i_k+N_k}(t)$}
\psfrag{qiNp1}[l]{$q_{i_k+N_k+1}(t)$}
\psfrag{chi=1}[l]{$\chi(t,\cdot)=1$}
\centerline{\resizebox{0.50\textwidth}{!}{
{\includegraphics{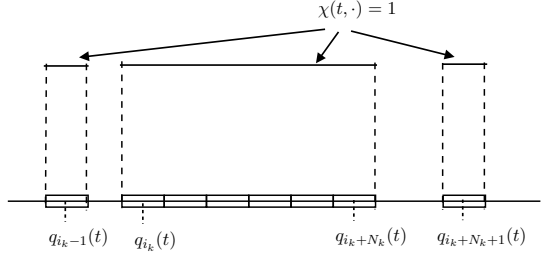}}
}}
\caption{Solid phase function at time $t$: $t\longrightarrow\chi(t,.)$.}
\label{fig:chit}
\end{figure}

\begin{figure}[hbtp]
\psfragscanon
\psfrag{qim1}[l]{$q_{i_k-1}(t)$}
\psfrag{qi}[l]{$q_{i_k}(t)$}
\psfrag{qiN}[l]{$q_{i_k+N_k}(t)$}
\psfrag{qiNp1}[l]{$q_{i_k+N_k+1}(t)$}
\psfrag{ui}[l]{$u_{i_k}(t)$}
\psfrag{uim1}[l]{$u_{i_k-1}(t)$}
\psfrag{uiNp1}[l]{$u_{i_k+N_k+1}(t)$}
\psfrag{uiN}[l]{$u_{i_k+N_k}(t)=u_{i_k}(t)$}
\centerline{\resizebox{0.80\textwidth}{!}{
{\includegraphics{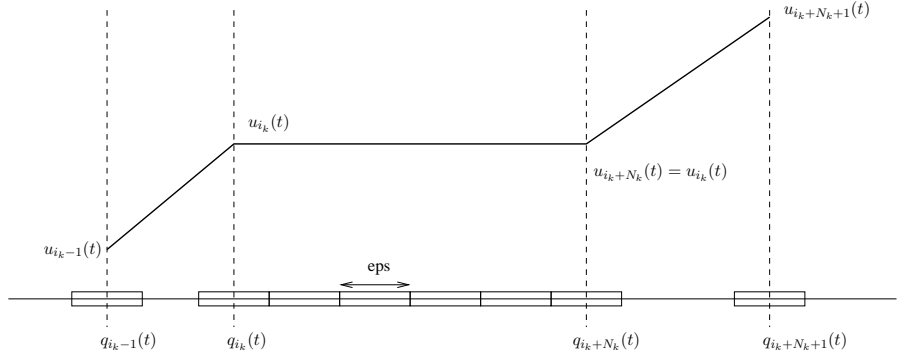}}
}}
\caption{Velocity function at time $t$: $t\longrightarrow u(t,.)$.}
\label{fig:ut}
\end{figure}

To show the convergence of the discrete model~(\ref{EDO}) to this
continuous model, we construct a micro-to-macro operator $\tQQ$. It is obtained from $\QQ$ in a similar way than $\tPP$ is defined from $\PP$  (see~(\ref{eq:udef})). Given $T>0$, an
$\eps$-feasible initial position vector $\qqinit\in\IR^{N-1}$ and a
force density $f\in L^1(]0,T[\times I)$, we define $\qq$ as $\QQ(\qqinit,\ff^\eps,\eps)$, where $\ff^\eps$ is given by
$$
f_i^\eps(t)=\frac{1}{2\eps}\int_{q_i(t)-\eps}^{q_i(t)+\eps} f(t,s)ds.
$$
We denote by $\uu=\dot\qq$ the velocities of the particles. Similarly
to what has been done in the previous section (see~(\ref{eq:udef})
and~(\ref{eq:chidef})), for each time $t$ we associate to vectors
$\qq(t)$ and $\uu(t)$ the functions $\chi(t,\cdot)$ and $u(t,\cdot)$
respectively constant by part and affine by part (see
figures~\ref{fig:chit} and~\ref{fig:ut}). We denote this mapping by  $\tilde{\mathcal
  Q}$~: $(\chi,u)=\tQQ(\bar \qq,f,\eps)$.

To state the convergence result it remains to say what is meant by
being a solution to the transport equation:
\begin{definition}
\label{def_eq_transport}
If $\rho_0\in\Linf$, we say that $$(\rho,u)\in\LinfLinf\times\LunLun$$
is a weak solution to problem
$$
\left|
\begin{array}{l}
\dsp \partial_t\rho +\partial_x (\rho u)=0,\vspace{6pt}\\
\dsp \rho(0,\cdot)=\rho_0,
\end{array}
\right.
$$
if, for all $\phi\in{\mathcal D}([0,T[\times I)$, 
\[
\int_0^T\int_I \rho(t,x)\partial_t\phi(t,x) dxdt 
+ \int_I\rho_0(x)\phi(0,x)dx
+ \int_0^T\int_I \rho(t,x)u(t,x)\partial_x\phi(t,x) dxdt
=0.
\]
\end{definition}

The convergence result in the non-stationary case is the following:

\begin{theorem}

Let $T>0$ be given and $f$ and $\rho_0$  be two measurable functions
on $]0,T[\times I$ and $I$ respectively. We suppose:
\begin{equation}
\label{H_f_int} f \in L^\infty(]0,T[,W^{1,\infty}(I))  
\cap
W^{1,\infty}(]0,T[,L^1(I))  .
\end{equation}
\begin{equation}
\label{H_rho} \rho_0 \in L^\infty(I) \hbox{ and } 0 \leq \rho_0(x) \leq 1\quad a.e. \hbox{ in }
I.
\end{equation}

Let $(\qqinit^\varepsilon)_\varepsilon$ be a sequence of
$\eps$-feasible vectors (see definition~\ref{definition:eps_feasible}), such that  
\begin{equation}
\label{rho0_cv}
\chi(\qqinit^\eps,\eps)\wsto \rho_0 \hbox{ in } \Linf
\end{equation}
where $\chi(\qqinit^\eps,\eps)$ is defined by~(\ref{eq:chidef}).
We denote by $(\chi^\eps,u^\eps)=\tilde{\mathcal Q}(\qqinit^\eps,f,\eps)$  the solution to the discrete problem.



There
exists $\rho \in \LinfLinf$ and
$u \in  \LinfHun\cap W^{1,\infty}(]0,T[,\Linf)$, such that, up to a subsequence,
\[
\begin{array}{l}
\dsp\chi^\varepsilon \wsto \rho \hbox{ in }\LinfLinf,
\vspace{6pt}\\
\dsp u^\varepsilon \wsto u \hbox{ in } \LinfHun,
\vspace{6pt}\\
\partial_t u^\eps \wsto \partial_t u \hbox{ in } \LinfLinf.
\end{array}
\]
Moreover, $\rho$ and $u$ verify
\begin{equation}
\label{rho_inf_1}
0 \leq \rho(t,x) \leq 1 \hbox{  for a.e. 
   } (x,t) \in ]0,T[\times I
\end{equation}
\begin{equation}
\label{eq_forcet}
-\partial_x\left(\frac{1}{1-\rho(t,\cdot)}\partial_xu(t,\cdot)\right)=\rho(t,\cdot) f(t,\cdot)
\hbox{    in the sense of
Prop.~\ref{proposition:def_sol_Kinfini}}
,
 \hbox{  for a.e. 
   } t \in ]0,T[
\end{equation}
\begin{equation}
\label{eq_transport}
\partial_t \rho + \partial_x(\rho u)=0
\virg
\rho(0,\cdot)=\rho_0
\hbox{ in the sense of Def.~\ref{def_eq_transport}}.
\end{equation}
\end{theorem}

\begin{proof}
The proof of this theorem is similar to the one of
theorem~\ref{thm:stat}. The main difference is that now, we need to control the time regularity of
the functions involved in order to pass to the limit when $\eps$ goes
to zero, and to check that the transport equation is verified. The
sketch of the proof is therefore similar to the one of
theorem~\ref{thm:stat}. Some of the lemmas involved here are even direct
consequences of the computations made to prove their stationnary
counterpart. These lemmas will be noticed in the
following. For readability reasons, the proofs of the remaining lemmas
will be postponed to the end of the section. \vspace{6pt}

To begin, note that hypothesis~(\ref{H_f_int}) for $f$,
together with proposition~\ref{prop:Q}, gives the existence of a
unique global solution to the discrete problem. \vspace{6pt}

The fact that $|\chi^\varepsilon|$ is bounded by $1$ immediatly gives
\begin{lemma}
\label{chi_cvt}
\[\exists\rho\in \LinfLinf \hbox{ such that } 
\chi^\varepsilon \wsto \rho \hbox{ in } \LinfLinf\]
\end{lemma}
Moreover, since $0\leq \chi^\eps\leq 1$, we
obtain~(\ref{rho_inf_1}).

Similarly to what has been done for the stationary case (see~(\ref{eq_rho_eps})
and~(\ref{eq_w_eps})), we define for each time
$t$, $\rho^\eps(t,\cdot)$ and $w^\eps(t,\cdot)$ (see figures~\ref{fig:rho_epst} and~\ref{fig:w_epst})
\begin{figure}[hbtp]
\psfragscanon
\psfrag{qim1}[l]{$q_{i_k-1}^\varepsilon(t)$}
\psfrag{qi}[l]{$q_{i_k}^\varepsilon(t)$}
\psfrag{qiN}[l]{$q_{i_k+N_k}^\varepsilon(t)$}
\psfrag{qiNp1}[l]{$q_{i_k+N_k+1}^\varepsilon(t)$}
\psfrag{rhoi}[l]{$\rho_{i_k}^\varepsilon(t)$}
\psfrag{rho1}[l]{$\rho^\varepsilon(t)=1$}
\psfrag{rhoiNp1}[l]{$\rho_{i_k+N_k+1}^\varepsilon(t)$}
\centerline{\resizebox{0.80\textwidth}{!}{
{\includegraphics{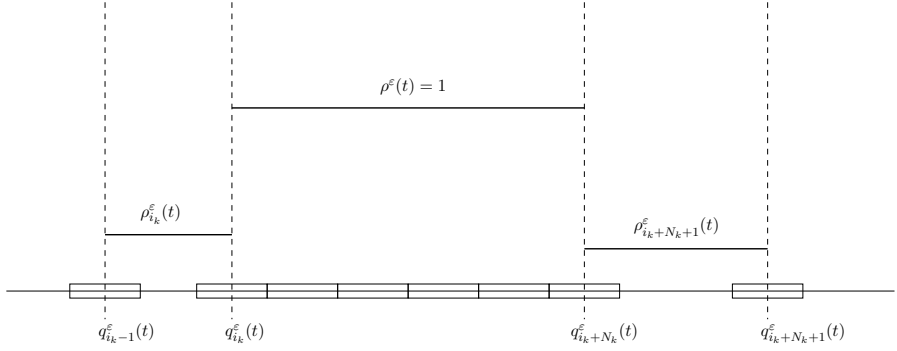}}
}}
\caption{Density function: $t\rightarrow\rho^\varepsilon(t,.)$.}
\label{fig:rho_epst}
\end{figure}
\begin{figure}[hbtp]
\psfragscanon
\psfrag{a}[l]{$0=q_{0}^\varepsilon$}
\psfrag{q1}[l]{$q_{1}^\varepsilon(t)$}
\psfrag{qim1}[l]{$q_{i_k-1}^\varepsilon(t)$}
\psfrag{qi}[l]{$q_{i_k}^\varepsilon(t)$}
\psfrag{qiN}[l]{$q_{i_k+N_k}^\varepsilon(t)$}
\psfrag{qiNp1}[l]{$q_{i_k+N_k+1}^\varepsilon(t)$}
\psfrag{qNm1}[l]{$q_{N^\eps-1}^\varepsilon(t)$}
\psfrag{b}[l]{$1=q_{N^\eps}^\varepsilon$}
\psfrag{w1}[l]{$w_{1}^\varepsilon(t)$}
\psfrag{w2}[l]{$w_{2}^\varepsilon(t)$}
\psfrag{wim1}[l]{$w_{i_k-1}^\varepsilon(t)$}
\psfrag{wi}[l]{$w_{i_k}^\varepsilon(t)$}
\psfrag{wiNp1}[l]{$w_{i_k+N_k+1}^\varepsilon(t)$}
\psfrag{Bip1}[l]{$w_{i_k+1}^\varepsilon(t)=\beta_{i_k+1}^\varepsilon(t)$}
\psfrag{BiN}[l]{$w_{i_k+N_k}^\varepsilon(t)=\beta_{i_k+N_k}^\varepsilon(t)$}
\psfrag{wNm1}[l]{$w_{N^\eps-1}^\varepsilon(t)$}
\psfrag{wN}[l]{$w_{N^\eps}^\varepsilon(t)$}
\centerline{\resizebox{0.95\textwidth}{!}{
{\includegraphics{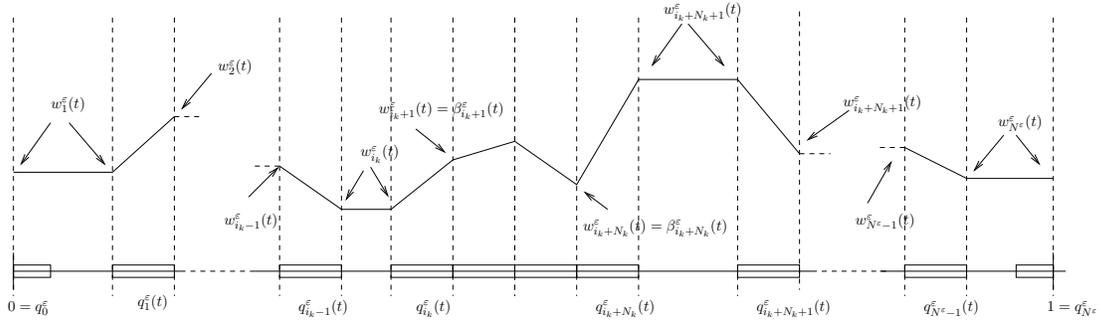}}
}}
\caption{$t\rightarrow w^\varepsilon(t,.)$.}
\label{fig:w_epst}
\end{figure}

From the computations carried out in the proof of
lemma~\ref{lem_u_eps} we get that~(\ref{eq_forcet}) is true at the
discrete level:
\begin{lemma}
\label{lem_u_epst} 

 For any $t\in ]0,T[$ and $\varepsilon>0$, 
\begin{equation}
\label{eq_force_epst}
-\partial_xw^\varepsilon(t,\cdot)=f^\varepsilon(t,\cdot),
\end{equation}
where $f^\varepsilon(t,\cdot) ={\dsp
\sum_{i=1}^{N^\eps-1}}f_i^\varepsilon(t)\ONE_{]q_i^\varepsilon(t)-\varepsilon,q_i^\varepsilon(t)+\varepsilon[}$.
\end{lemma}

The following lemma shows that the transport equation~(\ref{eq_transport}) is also true at
the discrete level:
\begin{lemma}
\label{lem_eq_transport} 

For any $\varepsilon>0$,
\begin{equation}
\label{eq_transport_eps}
\partial_t \rho^\varepsilon + \partial_x \left( \rho^\varepsilon
u^\varepsilon \right) =0,
\end{equation}
in the sense of definition~\ref{def_eq_transport}, with initial condition $\rho^\eps(0,\cdot)$.
\end{lemma}
\begin{remark} This result follows from the definition of  $\rho^\varepsilon$
and $u^\varepsilon$ with respect to $\qq^\varepsilon$.
\end{remark}

The proofs of lemmas~\ref{rho_cv},~\ref{f_cv} and~\ref{u_eps_bornee} extend straightforwardly to the non-stationary case, which yields
\begin{equation}
\label{eq:rho_cvt}
\rho^\varepsilon \wsto \rho\hbox {  in }\LinfLinf  \hbox { and } 
\rho^\varepsilon(0,.) \wsto \rho_0\hbox { in } \Linf
\end{equation}
$$
f^\varepsilon\wsto \rho f\hbox { in } \LinfLinf,
$$

\begin{equation}
\label{eq:u_eps_borneet}
(u^\varepsilon)_\varepsilon\hbox {  is bounded in }\LinfHun.
\end{equation}
From~(\ref{eq:u_eps_borneet}), we get the existence of a subsequence,
still denoted by $(u^\eps)_\eps$ such that, when $\eps$ goes to zero,
\begin{equation}
\label{u_eps_cv_faiblet} u^\varepsilon \wsto u \hbox{ in }
\LinfHun ,
\end{equation}
\begin{equation}
\label{dxu_eps_cv_faiblet} \partial_x u^\varepsilon \wsto \partial_x u \hbox{ dans }
\LinfLdeux .
\end{equation}

Similarly to the stationary case, in order to pass to the limit in the
left-hand side of~(\ref{eq_forcet}), we have to study the convergence
of $(w^\eps)_\eps$. The same computations as in the proof of
lemma~\ref{w_eps_bornee} give:
\begin{lemma}
\label{w_eps_borneet} 
$(w^\varepsilon)_\varepsilon$ and $(\partial_x
w^\varepsilon)_\varepsilon$ are bounded in $\LinfLinf$.
\end{lemma}

To obtain a strong convergence for $w^\eps$, we need to control its
time derivative, which will follow from the next lemma.
\begin{lemma}
\label{dt_w_eps_borneet} 
$(\partial_t w^\varepsilon)_\varepsilon$ is bounded $\LinfLinf$.
\end{lemma}
From lemmas~\ref{w_eps_borneet} and~\ref{dt_w_eps_borneet}, together
with Ascoli theorem, we get the existence of $w\in C^0(]0,T[\times I)$
    and a subsequence still denoted $(w^\eps)_\eps$ such that
\begin{equation}
\label{w_eps_cvt} w^\varepsilon \longrightarrow w \hbox{ in }
{\mathcal C}^0(]0,T[\times I) \hbox { when } \varepsilon\rightarrow 0.
\end{equation}
From this result of strong convergence we can establish, as in the
previous section, the following lemma:
\begin{lemma}
\label{dx_u_eps_cvt}
$\partial_x u^\varepsilon \wsto (1-\rho)w$ in $\LinfLinf$.
\end{lemma}
Finally, similarily to the stationary case, we pass to the limit
    in~(\ref{eq_force_epst}) to obtain~(\ref{eq_forcet}) by using test
    functions of type $(x,t)\rightarrow v(x)\Phi(t)$ with
$v\in\Hunzero$ and $\Phi\in \mathcal{D}(]0,T[)$.

To finish the proof of the theorem, it remains to
check~(\ref{eq_transport}) by passing to the limit
in~(\ref{eq_transport_eps}). To do so a strong convergence of
$(u^\eps)_\eps$ is required, which will follow from the next lemma:
\begin{lemma} 
\label{dt_u_eps_bornee}
$(\partial_t u^\varepsilon)_\varepsilon$ is bounded in $\LinfLinf$.
\end{lemma}
Using this last lemma, together with  the fact that $(u^\varepsilon)_\varepsilon$ is bounded in $\LinfHun$,
we
deduce that $(u^\eps)_\eps$ is bounded in $\HunHun$ and consequently,
there exists a subsequence, still denoted by $(u^\eps)_\eps$, such that
\begin{equation}
\label{u_eps_cv_fort} u^\varepsilon \longrightarrow u \hbox{ in }
\LdeuxLdeux  \hbox { when } \varepsilon\rightarrow 0.
\end{equation}
We are now going to let $\eps$ go to zero
in~(\ref{eq_transport_eps}). Let $\Psi$ be in ${\mathcal
  D}([0,T[\times I)$,
from lemma~\ref{lem_eq_transport}, we have
\[ \dsp \int_0^T\int_I \rho^\varepsilon\partial_t\Psi dxdt 
+\int_I \rho^\eps(0,x)\Psi(0,x)dx + \int_0^T\int_I
\rho^\varepsilon u^\varepsilon\partial_x\Psi dxdt = 0.
\]
The convergence of the first two terms comes from~(\ref{eq:rho_cvt}). To study the last term, merely write
 \[
\begin{array}{l}
\dsp
\left|\int_0^T\int_I
\rho^\varepsilon u^\varepsilon\partial_x\Psi dxdt-\int_0^T\int_I \rho
u\partial_x\Psi dxdt\right|
\vspace{6pt}\\
\dsp \quad\quad
\leq\left|\int_0^T\int_I\rho^\varepsilon(u^\varepsilon-u)\partial_x\Psi dxdt\right|
+\left|\int_0^T\int_I(\rho^\varepsilon-\rho)u\partial_x\Psi dxdt\right|,
\end{array}
\]
and use~(\ref{u_eps_cv_fort}) together with~(\ref{eq:rho_cvt}) to
show that it converges to zero, which completes the proof of the theorem.
\end{proof}

{\bf 
{Proof of Lemma \ref{lem_eq_transport}}}: $\partial_t \rho^\varepsilon
+ \partial_x \left( \rho^\varepsilon u^\varepsilon \right) =0$.

Let $\Phi$ be given in ${\mathcal D}([0,T[\times I)$. We first compute the
    time-derivative term. By using the definition~(\ref{eq_rho_eps}) of $\rho^\eps$
   we obtain:
\[
\int_0^T\int_a^b\rho^\eps(t,x)\partial_t\phi(t,x)dxdt = \sum_{i=1}^{N^\eps}\int_0^T\int_a^b\rho_i^\eps(t)\partial_t\phi(t,x)\ONE_{[q^\eps_{i-1}(t),q^\eps_i(t)]}(x)dxdt
\]
which gives
\[
\begin{array}{l}
\dsp \int_0^T\int_a^b\rho^\eps(t,x)\partial_t\phi(t,x)dxdt =\dsp 
-\sum_{i=1}^{N^\eps}\int_0^T\int_a^b(\rho_i^\eps)'(t)\phi(t,x)\ONE_{[q^\eps_{i-1}(t),q^\eps_i(t)]}(x)dxdt
\vspace{6pt}\\
\dsp \quad \quad -\sum_{i=1}^{N^\eps}\int_a^b\rho_i^\eps(0)\phi(0,x)\ONE_{[q^\eps_{i-1}(0),q^\eps_i(0)]}(x)dx
\vspace{6pt}\\
\dsp \quad \quad -\sum_{i=1}^{N^\eps}\int_0^T \left(
\rho_i^\eps(t)u^\eps_i(t)\phi(t,q^\eps_i(t)) - \rho_i^\eps(t)u^\eps_{i-1}(t)\phi(t,q^\eps_{i-1}(t))\right)dt.
\end{array}
\]
By shifting the indexes in the last sum and computing
$(\rho_i^\eps)'$,
we obtain
\begin{equation} 
\label{transport_dt}
\begin{array}{l}
\dsp \int_0^T\int_a^b\rho^\eps(t,x)\partial_t\phi(t,x)dxdt = \vspace{6pt}\\
\dsp \quad \quad \sum_{i=1}^{N^\eps}\int_0^T\int_a^b\frac{2\eps(u^\eps_i(t)-u^\eps_{i-1}(t))}{(q^\eps_i(t)-q^\eps_{i-1}(t))^2}\phi(t,x)\ONE_{[q^\eps_{i-1}(t),q^\eps_i(t)]}(x)dxdt
\vspace{6pt} \\
\dsp \quad \quad  -\int_a^b\rho^\eps(0)\phi(0,x)dx
+\sum_{i=1}^{N^\eps}\int_0^T 
(\rho_{i+1}^\eps(t)-\rho_i^\eps(t))u^\eps_i(t)\phi(t,q^\eps_i(t))dt.
\end{array}
\end{equation}

To compute the $x$-derivative term, an integration by part on each
sub-interval $[q_{i-1}^\eps(t),q_{i}^\eps(t)]$ gives
\begin{equation}
\label{transport_dx}
\begin{array}{l}
\dsp\int_0^T\int_a^b\rho^\eps(t,x)u^\eps(t,x)\partial_x\phi(t,x)dxdt 
\vspace{6pt}\\
\dsp\quad\quad=
\sum_{i=1}^{N^\eps}\int_0^T 
(\rho_{i}^\eps(t)-\rho_{i+1}^\eps(t))u^\eps_i(t)\phi(t,q^\eps_i(t))dt
\vspace{6pt}\\
\dsp\quad\quad\quad\quad-\sum_{i=1}^{N^\eps}\int_0^T\int_a^b\frac{2\eps(u^\eps_i(t)-u^\eps_{i-1}(t))}{(q^\eps_i(t)-q^\eps_{i-1}(t))^2}\phi(t,x)\ONE_{[q^\eps_{i-1}(t),q^\eps_i(t)]}(x)dxdt.
\end{array}
\end{equation}

The proof of the lemma is completed by summing up~(\ref{transport_dt}) and~(\ref{transport_dx}).
\qed

{\bf 
{Proof of Lemma \ref{dt_w_eps_borneet}}}: $(\partial_t
w^\varepsilon)_\varepsilon$ is bounded in $\LinfLinf$.

Using the definition~(\ref{eq_w_eps}) of $w^\eps$ (we recall that $w_i^\eps$ now depends on $t$), a simple
computation gives:
\[
\partial_t w^\eps(x,t)=w^\eps_{\hbox{\tiny der,aff}}(x,t)-w^\eps_{\hbox{\tiny der,cst}}(x,t),
\]
where, for each time $t$, $w^\eps_{\hbox{\tiny der,cst}}(t,\cdot)$ is the
piecewise constant function such that
\[
w^\eps_{\hbox{\tiny der,cst}}(x,t)=w^\eps_{\hbox{\tiny cst},i}(t)\stackrel{\hbox{\tiny def}}{=}
\frac{w^\eps_{i+1}(t)-w^\eps_i(t)}{2\eps}u^\eps_i(t)\quad \text{ on
}\quad ]q^\eps_i(t)-\eps,q^\eps_i(t)+\eps[,
\]
and $w^\eps_{\hbox{\tiny der,aff}}(t,\cdot)$ is the piecewise affine
function such that
\[
\begin{array}{ll}
\dsp w^\eps_{\hbox{\tiny der,aff}}(q^\eps_i(t)-\eps,t)=(w^\eps_i)'(t),\vspace{6pt}\\
\dsp w^\eps_{\hbox{\tiny der,aff}}(q^\eps_i(t)+\eps,t)=(w^\eps_{i+1})'(t).
\end{array}
\]
Therefore, to show that $(\partial_t w^\varepsilon)_\varepsilon$ is
bounded in $\LinfLinf$, it suffices to prove that the sets
$\{w^\eps_{\hbox{\tiny cst},i}(t)\}_i$, and $\{(w^\eps_i)'(t)\}_i$ are bounded
independently of $\eps$ and $t$.

By continuous injection of $\Hunzero$ in ${\mathcal C}(I)$ and
the boundedness of $(u^\varepsilon)_\varepsilon$ in $\LinfHun$,
we get that $(u^\eps)_\eps$ is bounded in
$\LinfLinf$ (we denote this bound by $M$). This, combined with the
fact that
$|w^\eps_{i+1}(t)-w^\eps_i(t)|=2\eps|f_i^\eps(t)|$ gives
$$
w^\eps_{\hbox{\tiny cst},i}(t)\leq\|u^\eps\|_\infty|f_i^\eps(t)|\leq M\|f\|_\infty,
$$
which gives a bound for $\{w^\eps_{\hbox{\tiny cst},i}(t)\}_i$.

To end the proof, it remains to show that $\{(w^\eps_i)'(t)\}_i$ is
bounded. As for the proof of lemma~\ref{w_eps_bornee}, it
suffices to show that 
$$\dsp\{(w^\eps_i)'(t)\}_{i, \text{ such that } \rho_i^\eps<1}$$ 
is bounded. From expression~(\ref{wi_eps}) of $w^\eps_i$, together
with the bound on $u^\eps$ and hypothesis~(\ref{H_f_int}) on $f$ we get
\begin{equation}
\label{borne_wiprim}
|(w^\eps_i)'(t)|
\leq M \|f\|_\LinfLun + MC(b-a) + \|\partial_t
f\|_\LinfLun,
\end{equation}
where $C$ is the Lipschitz constant of $f$, which completes the proof
of the lemma.
\qed

{\bf 
{Proof of Lemma \ref{dt_u_eps_bornee}}}: $(\partial_t u^\varepsilon)_\varepsilon$ is bounded in $\LinfLinf$.

As for the previous lemma, we compute $\partial_t u^\eps$ :
\[
\partial_t u^\eps = u^\eps_{\hbox{\tiny der,aff},1} - u^\eps_{\hbox{\tiny der,aff},2} - u^\eps_{\hbox{\tiny der,cste}},
\]
where, for each time $t$, $u^\eps_{\hbox{\tiny der,aff},1}(t,\cdot)$ is the
piecewise affine function such that
\[
u^\eps_{\hbox{\tiny der,aff},1}(t,q_i^\eps(t))=(u^\eps_{i})'(t),
\]
$u^\eps_{\hbox{\tiny der,aff},2}(t,\cdot)$ is the
piecewise affine function such that
\[
u^\eps_{\hbox{\tiny der,aff},2}(t,q_i^\eps(t))=\left(\frac{u^\eps_{i+1}(t)-u^\eps_{i}(t)}{q^\eps_{i+1}(t)-q^\eps_{i}(t)}\right)^2,
\]
and 
$u^\eps_{\hbox{\tiny der,cst}}(t,\cdot)$ is the
piecewise constant function such that
\[
u^\eps_{\hbox{\tiny der,cst}}(t,x)=\frac{u^\eps_{i+1}(t)-u^\eps_{i}(t)}{q^\eps_{i+1}(t)-q^\eps_{i}(t)}u^\eps_{i}(t)\quad \text{ on
}\quad ]q^\eps_i(t),q^\eps_{i+1}(t)[.
\]

Since $u^\eps$ is bounded (see proof of
lemma~\ref{dt_w_eps_borneet}), the fact that $u^\eps_{\hbox{\tiny der,aff},2}$ and
$u^\eps_{\hbox{\tiny der,cst}}$ are bounded comes from 
\begin{equation}
\label{borne2}
\left|\frac{u^\eps_{i+1}(t)-u^\eps_i(t)}{q^\eps_{i+1}(t)-q^\eps_i(t)}\right|=
(1-\rho^\eps_i(t))\left|w^\eps_i(t)\right|\leq\left|w^\eps_i(t)\right|\leq\|w^\eps\|_\LinfLinf,
\end{equation}
together with lemma~\ref{w_eps_borneet}.

It remains to prove that $u^\eps_{\hbox{\tiny der,aff},1}$ is bounded. Since, for
each $t$ we have $u^\eps_{\hbox{\tiny der,aff},1}(t,0)=0$, it suffices to show that
$\partial_x u^\eps_{\hbox{\tiny der,aff},1}$ is bounded. This in turn will come
provided we prove
$$\frac{(u^\eps_{i})'(t)-(u^\eps_{i-1})'(t)}{q^\eps_i(t)-q^\eps_{i-1}(t)}$$
is bounded independently of $t$, $\eps$ and $i$. To do so, we compute
the time-derivative of both side of the following equality
\[
\frac{u^\eps_{i}(t)-u^\eps_{i-1}(t)}{q^\eps_i(t)-q^\eps_{i-1}(t)}=(1-\rho^\eps_i(t))w^\eps_i(t)
\]
to obtain
\[
\left|\frac{(u^\eps_{i})'(t)-(u^\eps_{i-1})'(t)}{q^\eps_i(t)-q^\eps_{i-1}(t)}\right|\leq
|w^\eps_i(t)|^2+|(w^\eps_i)'(t)|+ \left(\frac{u^\eps_{i}(t)-u^\eps_{i-1}(t)}{q^\eps_i(t)-q^\eps_{i-1}(t)}\right)^2.
\]
Finally, from lemma~\ref{w_eps_borneet} and
inequalities~(\ref{borne_wiprim}) and~(\ref{borne2}), this last quantity is
bounded, which completes the proof of the lemma.
\qed

\section{Extensions}

A first extension to the approach we presented here consists in integrating inertial effects for the particles. Equation~(\ref{eq:sys}), which expresses instantaneous force balance, is replaced by Newton's  law
$$
\frac {d\uu}{dt} =  -A(\qq)\uu+ \ff 
$$
Having the number of masses go to infinity, it is natural to expect some quite of 1d pressureless Navier-Stokes equation:
$$
\partial_t (\rho u)+ \partial_t (\rho u^2)
- \partial_x \left( \frac{1
}{1-\rho}\partial_x u \right) = \rho f,
$$
coupled with the transport equation
$$
\partial_t \rho +\partial_x(\rho u) = 0.
$$
Now another natural question is the following: considering that fluid viscosity goes to $0$ (so that  the effective viscosity in the previous model is $\eps / (1-\rho)$, with $\eps\rightarrow 0$), what kind of limit model can be expected ? 
This very question has been addressed in~\cite{maurygluey} for  the case of a single particle against a wall. As $\eps$ goes to $0$, viscous effects are likely to disappear on unsaturated zones (where $\rho<1$). Besides, as
 lubrication forces (even for a vanishing viscosity) prevent contact in finite time, a model of the following type could be expected 
$$
\left |
\begin{array}{rcr}
\partial_t \rho  + \partial_x (\rho u)&=& 0\vspace{4pt}\\
\partial_t (\rho u)  + \partial_x (\rho u^2) + \partial _x p
&=& f\vspace{4pt}\\
\partial_t \gamma  + \partial_x (\gamma u) 
&=& - p \vspace{4pt}\\
\gamma \leq 0 \virg 
\rho \leq 1 \virg \gamma(1-\rho)& =& 0
\end{array}
\right .
$$
where $p$ is an unknown pressure-like field which may take negative values (which distin\-guishes this  gluey model from 
the standard pressureless gas situation), and $\gamma$ is a field which keeps track of the constraint history  experienced 
by a fluid  particle (see again~\cite{maurygluey} for more details on the gluey particle model). 

Extension of  this approach to higher dimensions is delicate. Firstly, even with strong assumptions on the shape of particles (discs or spheres), the notion of maximal  density is fuzzy. It may depend on the local arrangements of grains, and this structure at the microscopic level has to be described in some way at the macroscopic scale. 
Besides, for the same reason of local  geometrical complexity, the effects of lubrication forces cannot be expected to be described by a simple scalar (equivalent viscosity), but by a modification  of the stress tensor.  An asymptotic analysis is proposed in~\cite{berlyand,berlyand2} for the static problem, with a precise estimation of both shear and extensional viscosities, 
 but the problem of homogenization of the evolution problem is still widely open.


\begin{thebibliography}{99}

\bibitem{alibert}
J.J. Alibert, P. Seppecher \textit{Closure of the set of diffusion functionals - the one dimensional case}, Potential Analysis, 2008, no 28 pp 335-356. 

\bibitem{batchelor}
G.K. Batchelor, J.T. Green, 
J. Fluid Mech. {\bf 56}, p. 375, 1972.

\bibitem{bossis}
G. Bossis, J.F. Brady, J. Fluid Mech {\bf 155}, 1985.

\bibitem{briane1}
M. Briane, J. Casado-Diaz, \textit{Two-dimensional div-curl results. Application to the lack of nonlocal effects in homogenization}, Com. Part. Diff. Equa., 32 (2007), 935-969

\bibitem{briane2}
M. Briane, J. Casado-Dìaz, \textit{Asymptotic behaviour of equicoercive diffusion energies in two dimension}, Calc. Var. Part. Diff. Equa., 29  (4) (2007), 455-479. 

\bibitem{KK} Kim \& Karrila \textit{Microhydrodynamics : Principles
and Selected Applications}, Butterworth-Heinemann, Boston, 1991.

\bibitem{berlyand} L. Berlyand \& L. Borcea \& A. Panchenko,
\textit{Network approximation for effective viscosity of
concentred suspensions with complex geometry}, SIAM J. Math. Anal.
       {\bf 36} (2005), no5, 1580--1628.
       
  \bibitem{berlyand2}
L. Berlyan, A. Panchenko,  \textit{Strong and weak blow up of the viscous dissipation rates for concentrated suspensions}, Journal of Fluid Mechanics, v. 578, pp. 1-34 (2007).

\bibitem{brenner} Brenner,
\textit{Dissipation of Energy due to Solid Particles Suspended in a
  Viscous Liquid}, Physics of Fluid, 1958, Volume 1, Issue 4, pp 338-346.

\bibitem{einstein} A. Einstein, \textit{Eine neue Bestimmung der Molek\"uledimensionen}, {Ann. Phys. Leipzig} 1906, {\bf 19}, 289--306.

\bibitem{frankel} Frankel \& Akrivos,
\textit{On the Viscosity of a Concentrated Suspension of Solid
  Spheres}, Chemical Engineering Science, 1967, Vol. 22 pp. 847-853.

\bibitem{glo}
 R. Glowinski, T.W. Pan, T.I. Hesla, D.D. Joseph, J.  Périaux, {\em A fictitious domain approach to the direct numerical simulation of incompressible viscous flow past moving rigid bodies: application to particulate flow} J. Comput. Phys. 169 (2001), no. 2, 363--426.

\bibitem{landau}L. Landau, E. Lifschitz, \textit{Physique th\'eorique, m\'ecanique des fluides}, Ellipses.

\bibitem{mauryapp}
A. Lefebvre, B. Maury, {\em Apparent viscosity of a mixture of a Newtonian fluid and interacting particles},
Fluid-solid interactions: modeling, simulation, bio-mechanical applications, Comptes Rendus Mécanique, Volume 333, Issue 12, December 2005, p.p. 923-933.

\bibitem{maury} B.  Maury, \textit{A Many-Body Lubrication Model},  C. R. Acad. Sci. Paris, t. 325, S\'erie I, pp. 1053-1058, 1997.

\bibitem{maurygluey}
B. Maury, A gluey particle model, ESAIM Proceedings, July 2007, Vol.18, 133-142
Jean-Fr\'ed\'eric Gerbeau \& St\'ephane Labb\'e, Editors.

\bibitem{snabre} P.Mills, P. Snabre, \textit{Rheology and Strucuture of Concentrated Suspensions of Hard Spheres. Shear Induced Particle Migration},
J. Phys. II France, pp. 1597-1608, 1995.

\bibitem{murat} Murat, F.: H-convergence. S\'eminaire d'Analyse Fonctionnelle et Num\'erique, 1977-78, Universit\'e d'Alger. English translation : Murat F. \& Tartar L., H-convergence. Topics in the Mathematical Modelling of Composite Materials (Ed. Cherkaev, L. \& Kohn, R.V.). Progress in Nonlinear Differential Equations and their Applications, 31, Birka\:user, Boston, 21-43, 1998.

\bibitem{spagnolo} Spagnolo, S.: Sulla convergenza di soluzioni di equazioni paraboliche ed ellittiche. Ann. Sc. Norm. Super. Pisa Cl. Sci. (5) 22, 571-597 (1968).

\end{thebibliography}
\end{document}